\documentclass[journal]{IEEEtran}
\usepackage{amssymb}
\usepackage{epsfig}
\usepackage{graphicx}
\usepackage{colortbl}

\newtheorem{theorem}{Theorem}

\newtheorem{remark}{Remark}
\newtheorem{lemma}{Lemma}

\newtheorem{definition}{Definition}

\newtheorem{example}{Example}

\ifCLASSINFOpdf
\else
\fi
\begin{document}
    %
    \title{Exact Controllability of Discrete-Time Stochastic System with Multiplicative Noise}
    %
    %
    %

    \author{Juanjuan~Xu~ and ~Huanshui Zhang
        \thanks{*This work was supported by the National Natural Science Foundation of China under Grants 61821004, 62250056 and the Natural Science Foundation of Shandong Province (ZR2021ZD14, ZR2021JQ24), and High-level Talent Team Project of Qingdao West Coast New Area (RCTDJC-2019-05), Key Research and Development Program of Shandong Province (2020CXGC01208), and Science and Technology Project of Qingdao West Coast New Area (2019-32, 2020-20, 2020-1-4).}
        \thanks{J. Xu is with School of Control Science and Engineering, Shandong University, Jinan, Shandong, P.R. China 250061 (e-mail: juanjuanxu@sdu.edu.cn).}
        \thanks{H. Zhang is with College of Electrical Engineering and Automation, Shandong University of Science and Technology, Shandong, P.R. China 266590 (e-mail: hszhang@sdu.edu.cn).}
    }

    \maketitle

    \begin{abstract}

        This paper is concerned with the exact controllability of discrete-time stochastic system which is one of the basic problems of modern control theory. Though the exact controllability of continuous-time system governed by It\^{o} stochastic differential equations has been well studied in \emph{S. Peng, Progress in Natural Science, 1994}, the counterpart of the discrete-time case is still open due to the adaptiveness constraint of the controllers and the solvability challenging of stochastic difference equation with terminal value.
        The main contribution in this paper is to present both the Gramian matrix criterion and the Rank criterion for the exact controllability of discrete-time stochastic system. The novelty lies in the transformation of the forward stochastic difference equation into a novel backward one.

    \end{abstract}


    %
    \IEEEpeerreviewmaketitle

    \section{Introduction}

    \IEEEPARstart{C}{ontrollability} is one of the fundamental concepts in control theory and plays pivotal role in the analysis and design of dynamical control systems \cite{r1, r2}, which has wide applications in engineering \cite{r3}, machinery \cite{r4} and finance \cite{r5}. Roughly speaking, controllability generally means, that it is possible to steer dynamical control systems from an arbitrary initial state to a terminal state under the set of admissible controls \cite{r6}.
    Conceived by Kalman, the study of controllability was started at the beginning of 1960s \cite{r7}.
    Since then, plenty of essential progress has been made on the controllability problems.
    For instance, the fundamental work on controllability problem of classical linear systems has been proposed in \cite{r8}.
    The necessary and sufficient conditions were obtained in \cite{r9} for complete controllability of unperturbed impulse systems.
    \cite{r10} investigated the controllability and reachability criteria results for linear systems.
    In \cite{r11}, geometric criteria for controllability of switched systems was established.
    More details are referred to the survey paper \cite{r16,r17,r18,r19}.

    Noting that the aforementioned results are only related to deterministic systems.
    However, in practical application, there are many random factors in the system, such as fluctuating stock prices and physical systems subject to thermal fluctuations.
    Thus it is necessary to investigate the controllability of stochastic systems.
    \cite{r24} studied the stochastic controllability of linear stochastic systems and derived sufficient conditions in terms of Lyapunov method.
    Necessary and sufficient conditions were given in \cite{r} for different kinds of stochastic relative controllability of linear systems with delay in control.
    Stochastic approximate controllability was considered in \cite{r22}.
    \cite{r23} investigated partial approximate controllability of stochastic systems with controls only in the drift terms.
    In particular, necessary and sufficient conditions of exact controllability and exactly null controllability were given
    in \cite{r26} for stochastic systems governed by It\^{o} stochastic differential equation with the aid of the well-posedness of corresponding backward stochastic differential equations.
    Later, \cite{r27} studied the exact controllability of linear stochastic systems with bounded time-varying coefficients.
    \cite{r25} discussed the exact controllability of  linear stochastic systems with random coefficient.
    In addition, controllability of stochastic game-based control systems was studied in \cite{r31}. \cite{r28} investigated the
    controllability of forward and backward stochastic differential equations.
    It is noted that most results in the literature are focused on continuous-time case governed by It\^{o} stochastic differential equations.
    For the discrete-time stochastic system, only few results have been obtained for the stochastic controllability and stochastic observability, e.g., \cite{r20}, \cite{r21} and references therein. While little progress has been made on the exact controllability due to the adaptiveness constraint of the controller and the solvability challenging of stochastic difference equation with terminal value.

    In this paper, we study the exact controllability of discrete-time system governed by forward stochastic difference equations (FSDEs).
    The main contribution is to present both the Gramian matrix criterion and the Rank criterion for the exact controllability of the stochastic system. The key technique is to transform the FSDE into a novel backward stochastic difference equation (BSDE). Moreover, the results are applied to the cases with output measurements, input delay and state delay.

    The remainder of the paper is organized as follows. Section II presents the studied problem.
    Equivalent BSDE is shown in Section III. Exactly null controllability is discussed in Section IV.
    Exact controllability is studied in Section V. Exact controllability with input delay is shown in Section VI. Exact controllability with state delay is given in Section VII. Numerical examples are illustrated in Section VIII.
    Some concluding remarks are given in the last section.

    The following notations will be used throughout this paper: $R^n$
    denotes the set of $n$-dimensional vectors; $x'$ denotes the
    transpose of $x$; a symmetric matrix $M>0\ (\geq 0)$ means that $M$ is
    strictly positive-definite (positive semi-definite). $
    E[x(k)|\mathcal{F}(t-1)]$ denotes the conditional expectation with
    respect to the filtration $\mathcal{F}(t-1).$
    Denote $[0,N]=\{0, 1, \ldots, N\}$.


    \section{Problem Formulation}

    The stochastic system studied in this paper is in discrete time and the dynamics are governed by
    \begin{eqnarray}
        x(k+1)=[Ax(k)+Bu(k)]+w(k)[\bar{A}x(k)+\bar{B}u(k)],\label{c1}
    \end{eqnarray}
    where $x(k)\in R^n$ is the state, $u(k)\in R^{m}$ is the control input, $k$ is the discrete-time variable and takes values in $[0,T]$ where $T$ is finite or infinite. $A, B, \bar{A}, \bar{B}$ are constant matrices with compatible dimensions.
    The multiplicative noise $w(k)$ is a sequence of one-dimensional random variables defined on a complete probability space $\{\Omega,\mathcal{F},\mathcal{P}\}$,
    which is independent identical distribution (i.i.d)
    with zero mean and unit variance. $\mathcal{F}(k)$ is the the natural filtration
    generated by $w(k)$, i.e., $\mathcal{F}(k)=\sigma\{w(0), \ldots, w(k)\}$.
    Let $L^2(\Omega, R^s)$ represent the space of $R^s$-valued, square integrable random vectors on $\{\Omega,\mathcal{F},\mathcal{P}\}$,
    and $l_{\mathcal{F}}^2([0, T], R^s)$ the space of square summarabe stochastic processes $m(k)\in L^2(\Omega, R^s)$ which is $\mathcal{F}(k-1)$-adapted for any $k\in [0, T]$.

    Our aim is to derive the conditions for the exactly null controllability and the exact controllability of discrete-time multiplicative-noise system (\ref{c1}).

    To this end, we first present the definitions of the controllability as below.
    \begin{definition}\label{d1} (Exactly null controllability)
        System (\ref{c1}) is said to be exactly null controllable, if for any nonzero $x\in R^n$, there exist a positive integer $N$ and
        a sequence of controllers $\{u(k), k\in [0,N]\}$ which belongs to $l_{\mathcal{F}}^2([0, N], R^m)$ such that $x(0)=x$ and $x(N+1)=0$.
    \end{definition}

    \begin{definition}\label{d2} (Exact controllability)
        System (\ref{c1}) is said to be exactly controllable, if for any nonzero $x\in R^n$ and $\mathcal{F}(N)$-adapted $\xi\in\mathcal{S}$, there exist a positive integer $N$ and a sequence of controllers $\{u(k), k\in [0,N]\}$ which belongs to $l_{\mathcal{F}}^2([0, N], R^m)$ such that $x(0)=x$ and $x(N+1)=\xi$.
    \end{definition}

    \begin{remark}\label{rem3}
        In the deterministic case, i.e., $\bar{A}=0$ and $\bar{B}=0$, it follows that
        \begin{eqnarray}
            x(N+1)&=&A^{N+1}x(0)+\sum_{i=0}^NA^{N-i}Bu(i).\nonumber
        \end{eqnarray}
        This implies that if $A$ is nonsingular, the deterministic system (\ref{c1}) is controllable if and only if one of the conditions holds \cite{Control}:
        \begin{itemize}
            \item There exists a positive integer $N$ such that the matrix $Gram=\sum_{k=0}^NA^kBB'(A')^k$ is nonsingular.
            \item $Rank (\left[
            \begin{array}{ccc}
                B & \cdots & A^{n-1}B \\
            \end{array}
            \right])=n$.
        \end{itemize}
        However, it does not work for the stochastic case, i.e., $\bar{A}\neq0, \bar{B}\neq0$. In this case, it yields from (\ref{c1}) that
        \begin{eqnarray}
            x(N+1)&=&A(N)\cdots A(0)x(0)\nonumber\\
            &&+\sum_{i=0}^NA(N)\cdots A(i+1)B(i)u(i).\label{c5}
        \end{eqnarray}
        Noting that the white noises are involved in the matrices $A(N)\cdots A(0)$ and $A(N)\cdots A(i+1)B(i)$, the controllability condition like the deterministic case cannot be derived similarly.
    \end{remark}

    \begin{remark}\label{rem1}
        From (\ref{c5}), the state $x(k+1)$ for any $k$ is in the form of
        \begin{eqnarray}
            x(k+1)&=&[AA(k-1)\cdots A(0)x(0)\nonumber\\
            &&+A\sum_{i=0}^kA(k-1)\cdots A(i+1)B(i)u(i)]\nonumber\\
            &&+w(k)[\bar{A}A(k-1)\cdots A(0)x(0)\nonumber\\
            &&+\bar{A}\sum_{i=0}^kA(k-1)\cdots A(i+1)B(i)u(i)].\nonumber
        \end{eqnarray}
        It is noted that the state $x(k+1)$ depends linearly on $w(k)$. In other words, any $\mathcal{F}(k)$-adapted $R^n$-valued square integrable
        random vector can not be achieved for $x(k+1)$ like $w^2(k)\zeta$ where $\zeta\in R^n$.
        This motivates the constraint that the terminal state belongs to set $\mathcal{S}$ in Definition \ref{d2} which will be explicitly characterized in (\ref{c53}).
    \end{remark}

    \section{Equivalent BSDE}

    As analyzed in Remark \ref{rem3}, the traditional forward iteration loses efficacy for the controllability of the stochastic system (\ref{c1}).
    To this end, we assume that
    \begin{eqnarray}
        Rank(\bar{B})=n,\label{a3}
    \end{eqnarray}
    and equivalently transform the FSDE (\ref{c1}) into a kind of BSDEs which develops a new method to study the controllability of stochastic difference equation.

    Under (\ref{a3}), there exists an invertible matrix $M$ such that $\bar{B}M=\left[
    \begin{array}{cc}
        I & 0 \\
    \end{array}
    \right]
    $. Denote $BM=\left[
    \begin{array}{cc}
        L & F \\
    \end{array}
    \right]
    $ and make a transformation to the control input as
    \begin{eqnarray}u(k)=M\left[
        \begin{array}{c}
            q(k) \\
            v(k) \\
        \end{array}
        \right].\label{c24}
    \end{eqnarray}

    \begin{lemma}\label{lem1}
        Under (\ref{a3}) and the invertibility of the matrix $A-L\bar{A}$, FSDE (\ref{c1}) is equivalently reformulated as the following BSDE:
        \begin{eqnarray}
            x(k)=Cx(k+1)+\bar{C}z(k)+Dv(k)-Cw(k)z(k),\label{a5}
        \end{eqnarray}
        where $C\triangleq (A-L\bar{A})^{-1}, \bar{C}\triangleq-(A-L\bar{A})^{-1}L$ and $D\triangleq-(A-L\bar{A})^{-1}F$.

    \end{lemma}
    Proof. By using (\ref{c24}), (\ref{c1}) becomes
    \begin{eqnarray}
        x(k+1)&=&[Ax(k)+Lq(k)+Fv(k)]\nonumber\\
        &&+w(k)[\bar{A}x(k)+q(k)].\label{c6}
    \end{eqnarray}
    By letting $z(k)=\bar{A}x(k)+q(k)$, we have from (\ref{c6}) that
    \begin{eqnarray}
        x(k+1)&=&(A-L\bar{A})x(k)+Lz(k)+Fv(k)\nonumber\\
        &&+w(k)z(k).\label{c7}
    \end{eqnarray}
    By using the invertibility of the matrix $A-L\bar{A}$, (\ref{c7}) is equivalently rewritten as
    \begin{eqnarray}
        x(k)&=&(A-L\bar{A})^{-1}[x(k+1)-Lz(k)\nonumber\\
        &&-Fv(k)-w(k)z(k)].\nonumber
    \end{eqnarray}
    which is exactly (\ref{a5}). The proof is now completed.
    \hfill $\blacksquare$

    From Lemma \ref{lem1}, the exact controllability of system (\ref{c1}) is equivalent to that of system
    governed by BSDE (\ref{a5}). As comparison with the role of $z(k)$ as one part of control in system (\ref{c1}),
    $z(k)$ in BSDE (\ref{a5}) is one part of the solution $(x(k),z(k))$ iterated from the terminal state $x(N+1)$ as shown below.

    \begin{lemma}\label{lem2}
        Assume that (\ref{a5}) is solvable, then its solution is given by the following BSDE:
        \begin{eqnarray}
            x(k)&=&E\Big[C(k)x(k+1)+Dv(k)\Big|\mathcal{F}(k-1)\Big],\label{c8}\\
            z(k)&=&E[w(k)x(k+1)|\mathcal{F}(k-1)],\label{c2}
        \end{eqnarray}
        where $C(k)\triangleq C+\bar{C}w(k)$.

    \end{lemma}
    Proof. From (\ref{a5}), it follows that
    \begin{eqnarray}
        x(k)&=&E[Cx(k+1)+\bar{C}z(k)+Dv(k)|\mathcal{F}(k-1)]\nonumber\\
        &&+C\{x(k+1)-E[x(k+1)|\mathcal{F}(k-1)]\nonumber\\
        &&-w(k)z(k)\}.
    \end{eqnarray}
    By using the $\mathcal{F}(k-1)$-adaptedness of $x(k)$ and $z(k)$, the above equation can be reformulated as:
    \begin{eqnarray}
        x(k)&=&E[Cx(k+1)+\bar{C}z(k)+Dv(k)|\mathcal{F}(k-1)],\label{c3}\\
        0&=&x(k+1)-E[x(k+1)|\mathcal{F}(k-1)]\nonumber\\
        &&-w(k)z(k).\label{c4}
    \end{eqnarray}
    It yields from (\ref{c4}) that
    \begin{eqnarray}
        z(k)&=&E[w(k)x(k+1)|\mathcal{F}(k-1)],\nonumber
    \end{eqnarray}
    that is, (\ref{c2}) holds.
    By substituting (\ref{c2}) into (\ref{c3}), BSDE (\ref{c8}) follows.
    The proof is now completed.
    \hfill $\blacksquare$

    Based on Lemma \ref{lem2}, the solution $x(k)$ can be iteratively obtained from (\ref{c8}) with a given terminal value $x(N+1)$, and $z(k)$ can thus be derived from (\ref{c2}). However, we note that (\ref{c8}) is solvable for any $\mathcal{F}(N)$-adapted terminal value $x(N+1)$ while this is not the case for (\ref{a5}) where $x(k+1)$ must satisfy the form (\ref{c4}). To this end, we define the set $\mathcal{S}$ as
    \begin{eqnarray}
        \mathcal{S}&=&\{\xi|\xi\in L^2(\Omega, R^s)~ such~that~BSDE~(\ref{c52})~with~\nonumber\\
        &&~\tilde{x}_1(N+1)=\xi~is~solvable\}\label{c53}
    \end{eqnarray}
    where
    \begin{eqnarray}
        \tilde{x}_1(k)=C\tilde{x}_1(k+1)+\bar{C}\tilde{z_1}(k)-Cw(k)\tilde{z}_1(k).\label{c52}
    \end{eqnarray}

    With the definition of $\mathcal{S}$ in (\ref{c53}) and Definitions \ref{d1} and \ref{d2}, we have the equivalence among the exact controllability and the exactly null controllability for system (\ref{a5}).

    \begin{lemma}\label{lem4}
        Under (\ref{a3}) and the invertibility of the matrix $A-L\bar{A}$, system (\ref{a5}) is exactly controllable if and only if it is exactly null controllable.

    \end{lemma}
    Proof. On one hand, it is obvious that when $\tilde{x}_1(N+1)=0$, then $\tilde{x}_1(k)=0, \tilde{z}_1(k)=0$ for $k=0, 1, \cdots, N$. That is, $0\in \mathcal{S}$. Accordingly, the exact controllability implies the exactly null controllability for system (\ref{a5}).

    On the other hand, we prove that system (\ref{a5}) is exactly controllable if it is exactly null controllable.
    Firstly, from (\ref{c52}) with $\tilde{x}_1(N+1)=\xi\in \mathcal{S}$, we have
    \begin{eqnarray}
        \tilde{x}_1(0)=E[C(0)\cdots C(N)\xi].\label{c54}
    \end{eqnarray}

    Secondly, we introduce $(\tilde{x}_2(k), \tilde{z}_2(k))$ satisfying
    \begin{eqnarray}
        \tilde{x}_2(k)=C\tilde{x}_2(k+1)+\bar{C}\tilde{z}_2(k)+Dv(k)-Cw(k)\tilde{z}_2(k),\label{c30}
    \end{eqnarray}
    with $\tilde{x}_2(N+1)=0$. By using the exactly null controllability of (\ref{c1}), it follows that
    there exists a controller $v(k)$ such that $\tilde{x}_2(0)=x-E[C(0)\cdots C(N)\xi]$.

    Finally, by denoting $x(k)=\tilde{x}_1(k)+\tilde{x}_2(k), z(k)=\tilde{z}_1(k)+\tilde{z}_2(k)$ and using (\ref{c52}) and (\ref{c30}), we have that $x(k)$ obeys the BSDE (\ref{a5}) with $x(0)=\tilde{x}_1(0)+\tilde{x}_2(0)=x$ and $x(N+1)=\tilde{x}_1(N+1)+\tilde{x}_2(N+1)=\xi$.
    This indicates that the system (\ref{a5}) is exactly controllable.
    Thus, the equivalence between the exact controllability and the exactly null controllability of (\ref{a5}) follows.
    The proof is now completed. \hfill $\blacksquare$

    \section{Exact Controllability}

    Based on Lemma \ref{lem4}, we firstly study the exactly null controllability of system (\ref{a5}) by defining the following Gramian matrix:
    \begin{eqnarray}
        G_N&=&\sum_{i=0}^NE[C(0)C(1)\cdots C(i-1)D\nonumber\\
        &&\times D'C'(i-1)\cdots C'(0)].\label{c9}
    \end{eqnarray}

    \begin{theorem}\label{lem3}
        Under (\ref{a3}) and the invertibility of the matrix $A-L\bar{A}$, system (\ref{a5}) is exactly null controllable if and only if there exists a positive integer $N$ such that the Gramian matrix $G_N$ in (\ref{c9}) is invertible.
    \end{theorem}
    Proof. ``Sufficiency"
    By using the invertibility of the matrix $G_N$ in (\ref{c9}), we define the controller
    \begin{eqnarray}
        v(i)&=&D'C'(i-1)\cdots C'(0)G_N^{-1}x.\label{c10}
    \end{eqnarray}
    It is easily verified that (\ref{a5}) is solvable under the controller (\ref{c10}) with terminal condition $x(N+1)=0$.
    Accordingly, the solution can be iterated by using (\ref{c8}) as
    \begin{eqnarray}
        x(0)&=&E[C(0)x(1)+Dv(0)|\mathcal{F}(-1)]\nonumber\\
        &=&E[C(0)C(1)x(2)+C(0)Dv(1)+Dv(0)]\nonumber\\
        &=&E[C(0)C(1)\cdots C(N)x(N+1)\nonumber\\
        &&+C(0)C(1)\cdots C(N-1)Dv(N)+\cdots\nonumber\\
        &&+C(0)Dv(1)+Dv(0)]\nonumber\\
        &=&E[\sum_{i=0}^NC(0)C(1)\cdots C(i-1)Dv(i)].\label{c11}
    \end{eqnarray}
    By substituting (\ref{c10}) into (\ref{c11}), it follows that
    \begin{eqnarray}
        x(0)
        &=&E[\sum_{i=0}^NC(0)C(1)\cdots C(i-1)DD'\nonumber\\
        &&\times C'(i-1)\cdots C'(0) G_N^{-1}x]\nonumber\\
        &=&x.\nonumber
    \end{eqnarray}
    This implies that under $v(i)$ in (\ref{c10}), $x(k)$ in (\ref{a5}) satisfies $x(0)=x$ and $x(N+1)=0$,
    that is, the system (\ref{a5}) is exactly null controllable.

    ``Necessity" We now prove that if system (\ref{a5}) is exactly null controllable, then there exists $N>0$ such that the matrix $G_N$ in (\ref{c9}) is invertible. Otherwise, for any $N>0$, the matrix $G_N$ is singular, that is, there exists a nonzero vector $\eta$ such that $G_N \eta=0$. This gives that
    \begin{eqnarray}
        0&=&\eta'G_N \eta\nonumber\\
        &=&\eta'\sum_{i=0}^NE[C(0)C(1)\cdots C(i-1)D\nonumber\\
        &&\times D'C'(i-1)\cdots C'(0)]\eta.\nonumber
    \end{eqnarray}
    Since $\eta'C(0)C(1)\cdots C(i-1)DD'C'(i-1)\cdots C'(0)\eta\geq0$, it follows from the above equation that
    \begin{eqnarray}
        0&=&\eta'C(0)C(1)\cdots C(i-1)D, ~a.s.\label{c12}
    \end{eqnarray}
    By letting $x(0)=\eta$ and using the exact null controllability of system (\ref{a5}), there exists a control input $v(k)$ such that
    \begin{eqnarray}
        \eta&=&E[\sum_{i=0}^NC(0)C(1)\cdots C(i-1)Dv(i)].\nonumber
    \end{eqnarray}
    This further gives that
    \begin{eqnarray}
        \eta'\eta&=&E[\sum_{i=0}^N\eta'C(0)C(1)\cdots C(i-1)Dv(i)]\nonumber\\
        &=&0,\label{c13}
    \end{eqnarray}
    where (\ref{c12}) has been used in the derivation of the last equality. Recalling that $\eta$ is nonzero, this is a contradiction
    with (\ref{c13}).
    The proof is now completed. \hfill $\blacksquare$

    %

    \begin{remark}
        By using (\ref{c10}), one part $v(k)$ of the control input in (\ref{c24}) can be given by
        \begin{eqnarray}
            v(k)&=&D'C'(k-1)\cdots C'(0)G_N^{-1}x.\label{c26}
        \end{eqnarray}
        From $z(k)=\bar{A}x(k)+q(k)$ and $z(k)=E[w(k)x(k+1)|\mathcal{F}(k-1)]$,
        the other part $q(k)$ is given by
        \begin{eqnarray}
            q(k)&=&E[w(k)x(k+1)|\mathcal{F}(k-1)]-\bar{A}x(k).\label{c28}
        \end{eqnarray}
        In addition, from (\ref{c8}), one has
        \begin{eqnarray}\label{ceqR7}
            x(k)=E\big[\sum_{i=k}^{N}C(k)\cdots C(i-1)Dv(i)\big{|}\mathcal{F}(k-1)\big],
        \end{eqnarray}
        which implies that
        \begin{eqnarray}
            &&E[w(k)x(k+1){|}\mathcal{F}(k-1)]\nonumber\\
            &=&E\big[w(k)\sum_{i=k+1}^{N}C(k+1)\cdots C(i-1)Dv(i)\big{|}\mathcal{F}(k-1)\big].\nonumber
        \end{eqnarray}
        Together with (\ref{c28}), we have that
        \begin{eqnarray}\label{ceqR01}
            &&q(k)\nonumber\\
            \hspace{-3mm}&=&\hspace{-3mm}E\Big[w(k)\sum_{i=k+1}^{N}C(k+1)\cdots C(i-1)Dv(i)\Big{|}\mathcal{F}(k-1)\Big]\nonumber\\
            &&-\bar{A}x(k),
        \end{eqnarray}
        where $v(k)$ satisfies (\ref{c26}).
        To summarize, the control input (\ref{c24}) with $q(k)$ in (\ref{ceqR01}) and $v(k)$ in (\ref{c26})
        can ensure the exactly null controllability of system (\ref{c1}).
    \end{remark}

    Next, we give the Rank criterion for the exactly null controllability of (\ref{a5}) by defining
    \begin{eqnarray}
        R&=&\left[
        \begin{array}{cccccccc}
            D & CD & \bar{C}D & C^2D & \bar{C}^2D &  &  & \\
        \end{array}
        \right.\nonumber\\
        &&\left.
        \begin{array}{cccccccc}
            & &  & &  & C\bar{C}D & \bar{C}CD & \cdots \\
        \end{array}
        \right].\label{c48}
    \end{eqnarray}

    \begin{theorem}\label{thm1}
        Under (\ref{a3}) and the invertibility of the matrix $A-L\bar{A}$, system (\ref{a5}) is exactly null controllable if and only if
        \begin{eqnarray}
            Rank (R)=n.\nonumber
        \end{eqnarray}
    \end{theorem}
    \emph{Proof.} ``Necessity" Assume that system (\ref{a5}) is exactly null controllable, then there exist a positive integer $N$ and
    $\{v(k), k\in [0,N]\}$ such that $x(0)=x$ and $x(N+1)=0$. From Lemma \ref{lem2}, the solution of (\ref{a5}) under $v(k)$ can be given by (\ref{c8}) with the terminal state $x(N+1)=0$. We now prove $Rank(R)=n$ by contradiction. That is, assume that $Rank (R)<n$, then there exists a nonzero vector $\beta$ such
    that $\beta'R=0$.
    Define $Y(k)$ satisfying
    \begin{eqnarray}
        Y(k+1)&=&C'(k)Y(k),\label{c49}\\
        Y(0)&=&\beta.\label{c50}
    \end{eqnarray}
    It is obtained that
    \begin{eqnarray}
        Y(k)=C'(k-1)\cdots C'(0)\beta.\nonumber
    \end{eqnarray}
    Together with $\beta' R=0$ and (\ref{c50}), it follows that
    \begin{eqnarray}
        D' Y(k)=0.\label{c51}
    \end{eqnarray}
    In addition, by using (\ref{c49}) and the solution $x(k)$ given by (\ref{c8}), one has
    \begin{eqnarray}
        &&E[Y'(k)x(k)-Y'(k+1)x(k+1)]\nonumber\\
        &=&E\{Y'(k)E[C(k)x(k+1)+Dv(k)|\mathcal{F}(k-1)]\nonumber\\
        &&-Y'(k)C(k)x(k+1)\}\nonumber\\
        &=&E[Y'(k)Dv(k)]\nonumber\\
        &=&0,\nonumber
    \end{eqnarray}
    where (\ref{c51}) has been used in the derivation of the last equality.
    By taking summation with respect to $k$ from $0$ to $N$ on both sides of the above equation, it yields that
    \begin{eqnarray}
        E[Y'(0)x(0)-Y'(N+1)x(N+1)]=0.\nonumber
    \end{eqnarray}
    From $x(N+1)=0$ and $Y(0)=\beta$, it further gives that $\beta'x(0)=0$. This implies that the initial state $x(0)$ is not arbitrary for system (\ref{a5}) with $x(N+1)=0$ since there is a nonzero vector $\beta$ such that $\beta'x(0)=0$, which is a contradiction with the exactly null controllability of system (\ref{a5}). Thus, $Rank (R)=n$.

    ``Sufficiency" We now prove $Rank (R)=n$ implies the exactly null controllability of system (\ref{a5}). Firstly, we prove that if $Rank (R)=n$, then there exists a positive integer $N$ such that the Gramian matrix $G_N$ in (\ref{c9}) is invertible. Otherwise, for any positive integer $N$, the matrix $G_N$ is singular, that is, there exists a nonzero vector $\eta$ such that $G_N\eta=0$, i.e., $\eta'\sum_{i=0}^NE[C(0)C(1)\cdots C(i-1)DD'C'(i-1)\cdots C'(0)]\eta=0$ for any $N$. This implies that $\eta'E[C(0)C(1)\cdots C(i-1)DD'C'(i-1)\cdots C'(0)]\eta=0$ for any $i=0, 1, 2, \ldots,$ that is,
    \begin{eqnarray}
        0&=&\eta'DD'\eta,\nonumber\\
        0&=&\eta'[CDD'C+\bar{C}DD'\bar{C}']\eta\nonumber\\
        &=&\eta'\left[
        \begin{array}{cc}
            CD & \bar{C}D \\
        \end{array}
        \right]\left[
        \begin{array}{cc}
            CD & \bar{C}D \\
        \end{array}
        \right]'\eta,\nonumber\\
        0&=&\eta'[C^2DD'(C')^2+\bar{C}^2DD'(\bar{C}')^2+C\bar{C}DD'\bar{C}'C'\nonumber\\
        &&+\bar{C}CDD'C'\bar{C}']\eta\nonumber\\
        &=&\eta'\left[
        \begin{array}{cccc}
            C^2D & \bar{C}^2D & C\bar{C}D & \bar{C}CD \\
        \end{array}
        \right]\nonumber\\
        &&\times \left[
        \begin{array}{cccc}
            C^2D & \bar{C}^2D & C\bar{C}D & \bar{C}CD \\
        \end{array}
        \right]'\eta,\nonumber\\
        &\vdots&\nonumber
    \end{eqnarray}
    This gives that $\eta'R=0$, which is a contradiction with the fact that $Rank(R)=n$. Accordingly, there exists a positive integer $N$ such that the Gramian matrix $G_N$ in (\ref{c9}) is invertible. Combining with Theorem \ref{lem3}, system (\ref{a5}) is exactly null controllable.
    The proof is now completed. \hfill $\blacksquare$

    To summarize, we have the following equivalent conditions for the controllability of system (\ref{c1}).
    \begin{theorem}\label{thm2}
        Under (\ref{a3}) and the invertibility of the matrix $A-L\bar{A}$, the following results are equivalent:
        \begin{itemize}
            \item[1.] System (\ref{c1}) is exactly controllable.
            \item[2.] System (\ref{c1}) is exactly null controllable.
            \item[3.] There exists a positive integer $N$ such that the Gramian matrix $G_N$ in (\ref{c9}) is invertible.
            \item[4.] $Rank(R)=n$.
        \end{itemize}
    \end{theorem}
    Proof. The results follow directly from Lemma \ref{lem1}, Lemma \ref{lem4}, Theorem \ref{lem3} and Theorem \ref{thm1}. \hfill $\blacksquare$

    \section{Stochastic Systems with Output Measurements}

    In Section IV, the conditions for the exact controllability of system (\ref{c1}) have been obtained as concluded in Theorem \ref{thm2}.
    In this section, we will study the exact controllability of the system (\ref{c1}) with output measurements and apply the derived results to solve the exact controllability of system (\ref{c1}) with $Rank(\bar{B})<n$. In particular, the measurement of system (\ref{c1}) is given by
    \begin{eqnarray}
        y(k)=Hx(k),\label{c71}
    \end{eqnarray}
    where $H\in R^{l\times n}$ is of full row rank. In this case, the definition of partially exact controllability controllability for (\ref{c1}) with (\ref{c71}) is given as below.
    \begin{definition}\label{d3} ($H$-partially exactly null controllability)
        System (\ref{c1}) is said to be $H$-partially exactly null controllable, if for any nonzero $\bar{x}_1\in R^l$, there exist a positive integer $N$ and a sequence of controllers $\{u(k), k\in [0,N]\}$ which belongs to $l_{\mathcal{F}}^2([0, N], R^m)$ such that $Hx(0)=\bar{x}_1$ and $Hx(N+1)=0$.
    \end{definition}

    \begin{remark}
        The matrix $H$ is assumed to be of full row rank. This is consistent with the arbitrariness of the initial value $\bar{x}_1\in R^l$ by using $\bar{x}_1=Hx(0)$. Otherwise, if $rank (H)=r<l$, we can find a full-row-rank matrix $\bar{H}\in R^{r\times n}$ and study the $\bar{H}$-partially exactly null controllability. In fact, in this case, there exists an invertible matrix $T\in R^{l\times l}$ such that $TH=\left[
        \begin{array}{c}
            \bar{H} \\
            0 \\
        \end{array}
        \right]$.
        Then, $THx(0)=\left[
        \begin{array}{c}
            \bar{H}x(0)\\
            0 \\
        \end{array}
        \right]$ and $THx(N+1)=\left[
        \begin{array}{c}
            \bar{H}x(N+1)\\
            0 \\
        \end{array}
        \right]$.
        This implies that $\bar{H}x(N+1)=0$ is equivalent to $Hx(N+1)=0$, and $\bar{H}x(0)\in R^r$ is arbitrary when $x(0)$ can be arbitrarily chosen in $R^n$.
    \end{remark}

    Different from the exact controllability, only partial state $Hx(k)$ is concerned whose dynamics are given from Lemma \ref{lem1} as follows:
    \begin{eqnarray}
        Hx(k)&=&HCx(k+1)+H\bar{C}z(k)+HDv(k)\nonumber\\
        &&-HCw(k)z(k).\label{c56}
    \end{eqnarray}
    We now present the condition to ensure the $H$-partially exactly null controllability.
    \begin{theorem}\label{thm3}
        Under (\ref{a3}) and the invertibility of the matrix $A-L\bar{A}$, assume that there exist matrices $\mathcal{C}_1$ and $\bar{\mathcal{C}}_1$
        such that
        \begin{eqnarray}
            HC=\mathcal{C}_1H,~~H\bar{C}=\bar{\mathcal{C}}_1H,\label{c57}
        \end{eqnarray}
        then the following results are equivalent:
        \begin{itemize}
            \item[1.] System (\ref{c1}) is $H$-partially exactly null controllable.
            \item[2.] System (\ref{c1}) is $H$-partially exactly controllable, that is, for any nonzero $\bar{x}_1\in R^l$ and $\xi_1\in \mathcal{S}_1$, there exist a positive integer $N$ and a sequence of controllers $\{u(k), k\in [0,N]\}$ which belongs to $l_{\mathcal{F}}^2([0, N], R^m)$ such that $Hx(0)=\bar{x}_1$ and $Hx(N+1)=\xi_1$ where \begin{eqnarray}
                \mathcal{S}_1&=&\{\xi_1|\xi_1\in L^2(\Omega, R^s)~ such~that~BSDE\nonumber\\
                &&\bar{x}_1(k)=\mathcal{C}_1\bar{x}_1(k+1)+\bar{\mathcal{C}}_1\bar{z}_1(k)\nonumber\\
                &&-\mathcal{C}_1w(k)\bar{z}_1(k)~with~\bar{x}_1(N+1)=\xi_1\nonumber\\
                &&is~solvable\}.\nonumber
            \end{eqnarray}
            \item[3.] There exists a positive integer $N$ such that the Gramian matrix $\bar{G}_N$ is invertible where $\bar{G}_N$ is defined by \begin{eqnarray}
                \bar{G}_N&=&\sum_{i=0}^NE[\mathcal{C}_1(0)\mathcal{C}_1(1)\cdots \mathcal{C}_1(i-1)HD\nonumber\\
                &&\times D'H'\mathcal{C}_1'(i-1)\cdots \mathcal{C}_1'(0)],
            \end{eqnarray}
            where $\mathcal{C}_1(k)\triangleq \mathcal{C}_1+w(k)\bar{\mathcal{C}}_1$.
            \item[4.] $Rank (\bar{R})=l$ where $\bar{R}$ is defined by
            \begin{eqnarray}
                \bar{R}&=&\left[
                \begin{array}{cccccccc}
                    HD & \mathcal{C}_1HD & \bar{\mathcal{C}}_1HD & \mathcal{C}_1^2HD &  &&  & \\
                \end{array}
                \right.\nonumber\\
                &&\hspace{-12mm}\left.
                \begin{array}{cccccccc}
                    &  &  & & \bar{\mathcal{C}}_1^2HD&  \mathcal{C}_1\bar{\mathcal{C}}_1HD  & \bar{\mathcal{C}}_1\mathcal{C}_1HD & \cdots \\
                \end{array}
                \right].\nonumber
            \end{eqnarray}
        \end{itemize}
    \end{theorem}
    Proof. By using (\ref{c57}), (\ref{c56}) is reformulated as
    \begin{eqnarray}
        Hx(k)&=&\mathcal{C}_1Hx(k+1)+\bar{\mathcal{C}}_1Hz(k)+HDv(k)\nonumber\\
        &&-\mathcal{C}_1Hw(k)z(k).\nonumber
    \end{eqnarray}
    By letting $Hx(k)=\hat{x}_1(k)$ and $Hz(k)=\hat{z}_1(k)$, it follows that
    \begin{eqnarray}
        \hat{x}_1(k)&=&\mathcal{C}_1\hat{x}_1(k+1)+\bar{\mathcal{C}}_1\hat{z}_1(k)+HDv(k)\nonumber\\
        &&-\mathcal{C}_1w(k)\hat{z}_1(k).\nonumber
    \end{eqnarray}
    By using Theorem \ref{thm2}, the equivalence between Items 1, 2, 3 and 4 follows directly. \hfill $\blacksquare$

    We now apply the $H$-partially exactly controllability to study the system (\ref{c1}) with $Rank(\bar{B})<n$. Without loss of generality, let $\bar{B}=\left[
    \begin{array}{cc}
        I_r & 0 \\
        0 & 0 \\
    \end{array}
    \right]
    $ where $r<\min\{m,n\}$.
    In this case, by denoting $A=\left[
    \begin{array}{cc}
        A_{11} & A_{12} \\
        A_{21} & A_{22} \\
    \end{array}
    \right], B=\left[
    \begin{array}{cc}
        B_1 & B_2 \\
    \end{array}
    \right]
    , B_i=\left[
    \begin{array}{c}
        B_{1i} \\
        B_{2i} \\
    \end{array}
    \right], i=1,2, \bar{A}=\left[
    \begin{array}{cc}
        \bar{A}_{11} & \bar{A}_{12} \\
        \bar{A}_{21} & \bar{A}_{22} \\
    \end{array}
    \right], x(k)=\left[
    \begin{array}{c}
        x_1(k) \\
        x_2(k) \\
    \end{array}
    \right], u(k)=\left[
    \begin{array}{c}
        u_1(k) \\
        u_2(k) \\
    \end{array}
    \right]
    $, system (\ref{c1}) is reduced to
    \begin{eqnarray}
        x_1(k+1)&=&[A_{11}x_1(k)+A_{12}x_2(k)+B_{11}u_1(k)\nonumber\\
        &&+B_{12}u_2(k)]+w(k)[\bar{A}_{11}x_1(k)\nonumber\\
        &&+\bar{A}_{12}x_2(k)+u_1(k)],\label{c22}\\
        x_2(k+1)&=&[A_{21}x_1(k)+A_{22}x_2(k)+B_{21}u_1(k)\nonumber\\
        &&+B_{22}u_2(k)]+w(k)[\bar{A}_{21}x_1(k)\nonumber\\
        &&+\bar{A}_{22}x_2(k)].\label{c23}
    \end{eqnarray}
    For simplicity of discussions, we let $n-r=r$, $\bar{A}_{21}=I_{n-r}, \bar{A}_{22}=0$,
 and make the following denotations:
    \begin{eqnarray}
        &&\mathcal{A}_{11}=A_{11}-B_{11}\bar{A}_{11}, \mathcal{A}_{12}=A_{12}-B_{11}\bar{A}_{12},\nonumber\\
        &&\mathcal{A}_{21}=A_{21}-B_{21}\bar{A}_{11}, \mathcal{A}_{22}=A_{22}-B_{21}\bar{A}_{12},\nonumber\\
        &&\mathbb{A}=\left[
        \begin{array}{cc}
            \mathcal{A}_{11} & \mathcal{A}_{12} \\
            \mathcal{A}_{21} & \mathcal{A}_{22}  \\
        \end{array}
        \right]^{-1},\nonumber\\
        &&\mathbb{B}=-\mathbb{A}\left[
        \begin{array}{c}
            B_{11} \\
            B_{21} \\
        \end{array}
        \right],\mathbb{D}=-\mathbb{A}\left[
        \begin{array}{c}
            B_{12} \\
            B_{22} \\
        \end{array}
        \right].\nonumber\\
    \end{eqnarray}
    By letting \begin{eqnarray}
        u_1(k)=z_1(k)-\bar{A}_{11}x_1(k)-\bar{A}_{12}x_2(k),\label{c68}
    \end{eqnarray}
    (\ref{c22}) is rewritten as
    \begin{eqnarray}
        x_1(k+1)&=&[\mathcal{A}_{11}x_1(k)+\mathcal{A}_{12}x_2(k)+B_{11}z_1(k)\nonumber\\
        &&+B_{12}u_2(k)]+w(k)z_1(k),\label{c69}
    \end{eqnarray}
    and (\ref{c23}) is rewritten as
    \begin{eqnarray}
        x_2(k+1)&=&[\mathcal{A}_{21}x_1(k)+\mathcal{A}_{22}x_2(k)+B_{21}z_1(k)\nonumber\\
        &&+B_{22}u_2(k)]+w(k)x_1(k).\label{c70}
    \end{eqnarray}
    By using the definition of matrices $\mathbb{A}, \mathbb{B}$ and $\mathbb{D}$, we can reformulate (\ref{c69})-(\ref{c70}) as
    \begin{eqnarray}
        \left[
        \begin{array}{c}
            x_1(k) \\
            x_2(k) \\
        \end{array}
        \right]&=&\mathbb{A}\left[
        \begin{array}{c}
            x_1(k+1) \\
            x_2(k+1) \\
        \end{array}
        \right]+\mathbb{B}z_1(k)+\mathbb{D}u_2(k)\nonumber\\
        &&-w(k)\mathbb{A}\left[
        \begin{array}{c}
            z_1(k) \\
            x_1(k) \\
        \end{array}
        \right].
    \end{eqnarray}
    The controllability results are thus illustrated as follows.
    \begin{theorem}\label{lem5}
        If there exist $\mathbb{A}_1$ and $\mathbb{B}_1$
        such that
        \begin{eqnarray}
            \left[
            \begin{array}{cc}
                I_r & 0 \\
            \end{array}
            \right]
            \mathbb{A}=\mathbb{A}_1\left[
            \begin{array}{cc}
                I_r & 0 \\
            \end{array}
            \right],~~\left[
            \begin{array}{cc}
                I_r & 0 \\
            \end{array}
            \right]\mathbb{B}=\mathbb{B}_1\left[
            \begin{array}{cc}
                I_r & 0 \\
            \end{array}
            \right],\nonumber
        \end{eqnarray}
        then the following results are equivalence:
        \begin{itemize}
            \item[1.] System (\ref{c1}) is $\left[
            \begin{array}{cc}
                I_r & 0 \\
            \end{array}
            \right]$-partially exactly null controllable.
            \item[2.]  System (\ref{c1}) is $\left[
            \begin{array}{cc}
                I_r & 0 \\
            \end{array}
            \right]$-partially exactly controllable, that is, for any nonzero $\hat{x}\in R^r$ and $\hat{\xi}\in \hat{\mathcal{S}}$, there exist a positive integer $N$ and a sequence of controllers $\{u(k), k\in [0,N]\}$ which belongs to $l_{\mathcal{F}}^2([0, N], R^m)$ such that $\left[
            \begin{array}{cc}
                I_r & 0 \\
            \end{array}
            \right]x(0)=\hat{x}$ and $\left[
            \begin{array}{cc}
                I_r & 0 \\
            \end{array}
            \right]x(N+1)=\hat{\xi}$ where \begin{eqnarray}
                \hat{\mathcal{S}}&=&\{\hat{\xi}|\hat{\xi}\in L^2(\Omega, R^s)~ such~that~BSDE\nonumber\\
                &&\hat{x}(k)=\mathbb{A}_1\hat{x}(k+1)+\mathbb{B}_1\hat{z}(k)\nonumber\\
                &&-\mathbb{A}_1w(k)\hat{z}(k)~with~\hat{x}(N+1)=\hat{\xi}\nonumber\\
                &&is~solvable\}.\nonumber
            \end{eqnarray}

            \item[3.] There exists a positive integer $N$ such that the Gramian matrix $\hat{G}_N$ is invertible where $\hat{G}_N$ is defined by \begin{eqnarray}
                \hat{G}_N&=&\sum_{i=0}^NE[\mathbb{A}_1(0)\mathbb{A}_1(1)\cdots \mathbb{A}_1(i-1)\mathbb{D}_1\nonumber\\
                &&\times \mathbb{D}_1'\mathbb{A}_1'(i-1)\cdots \mathbb{A}_1'(0)],
            \end{eqnarray}
            where $\mathbb{A}_1(k)\triangleq \mathbb{A}_1+w(k)\mathbb{B}_1, \mathbb{D}_1\triangleq \left[
            \begin{array}{cc}
                I_r & 0 \\
            \end{array}
            \right]\mathbb{D}$.
            \item[4.] $Rank (\bar{R})=l$ where $\bar{R}$ is defined by
            \begin{eqnarray}
                \bar{R}&=&\left[
                \begin{array}{cccccccc}
                    \mathbb{D}_1 & \mathbb{A}_1\mathbb{D}_1 & \mathbb{B}_1\mathbb{D}_1 & \mathbb{A}_1^2\mathbb{D}_1 & \mathbb{B}_1^2\mathbb{D}_1 & \\
                \end{array}
                \right.\nonumber\\
                &&\left.
                \begin{array}{cccccccc}
                    &  &  &  &\mathbb{A}_1\mathbb{B}_1\mathbb{D}_1  &\mathbb{B}_1\mathbb{A}_1\mathbb{D}_1   & \cdots \\
                \end{array}
                \right].
            \end{eqnarray}
        \end{itemize}

    \end{theorem}
    Proof. The proof follows directly from Theorem \ref{thm3}. \hfill $\blacksquare$


    \section{Stochastic Systems with Input Delay}

    In this section, we consider the stochastic systems with input delay where the dynamics are governed by
    \begin{eqnarray}
        x(k+1)&=&[Ax(k)+B_1u_1(k-\tau)+Bu(k)]\nonumber\\
        &&+w(k)[\bar{A}x(k)+\bar{B}u(k)],\label{ct10}
    \end{eqnarray}
    where constant matrices $A, B, \bar{A}, \bar{B}$ have been defined in (\ref{c1}), $B_1$ is a constant matrix with compatible dimension, $\tau>0$ represents the input delay.

    By applying (\ref{a3}) and (\ref{c24}), (\ref{ct10}) becomes
    \begin{eqnarray}
        x(k+1)&=&[Ax(k)+B_1u_1(k-\tau)+Lq(k)+Fv(k)]\nonumber\\
        &&+w(k)[\bar{A}x(k)+q(k)]\nonumber\\
        &\triangleq&[(A-L\bar{A})x(k)+B_1u_1(k-\tau)\nonumber\\
        &&+Lz(k)+Fv(k)]+w(k)z(k),\nonumber
    \end{eqnarray}
    where $z(k)=\bar{A}x(k)+q(k).$

    \begin{lemma}
        Under (\ref{a3}) and the invertibility of the matrix $A-L\bar{A}$, system (\ref{ct10}) is equivalently rewritten as
        \begin{eqnarray}
            x(k)&=&Cx(k+1)+D_1u_1(k-\tau)+\bar{C}z(k)\nonumber\\
            &&+Dv(k)-Cw(k)z(k),\label{ct11}
        \end{eqnarray}
        where $C, \bar{C}, D$ have the same definitions as those in Lemma \ref{lem1} and $D_1\triangleq -(A-L\bar{A})^{-1}B_1$. In particular,
        the solution to (\ref{ct11}) is given by
        \begin{eqnarray}
            x(k)&=&E\Big[C(k)x(k+1)+D_1u(k-\tau)\nonumber\\
            &&+Dv(k)\Big|\mathcal{F}(k-1)\Big],\label{ct12}
        \end{eqnarray}
        and $z(k)$ satisfies (\ref{c2}).
    \end{lemma}
    Proof. The proof follows similarly to those in Lemma \ref{lem1} and \ref{lem2}. So we omit the details. \hfill $\blacksquare$

    By defining the Gramian matrix
    \begin{eqnarray}
        G^{\tau}_N&=&G_N+\sum_{i=0}^NE\Big[E[C(0)C(1)\cdots C(i-1)|\mathcal{F}(i-\tau-1)]\nonumber\\
        &&\hspace{-3mm}\times D_1D_1'E[C(0)C(1)\cdots C(i-1)|\mathcal{F}(i-\tau-1)]'\Big],\label{ct14}
    \end{eqnarray}
    the controllability condition for system (\ref{ct10}) can be obtained as shown below.

    \begin{theorem}\label{thm4}
        Under (\ref{a3}) and the invertibility of the matrix $A-L\bar{A}$, if there exists a positive integer $N$ such that the Gramian matrix $G^{\tau}_{N}$ in (\ref{ct14}) is invertible, then system (\ref{ct10}) with input delay is exactly null controllable, and is also exactly controllable, that is, for any nonzero $x\in R^n$ and $\xi\in \mathcal{S}$ defined by (\ref{c53}), there exist a positive integer $N$ and a sequence of controllers $\{u(k), k\in [0,N]\}$ which belongs to $l_{\mathcal{F}}^2([0, N], R^m)$ such that $x(0)=x$ and $x(N+1)=\xi$.
    \end{theorem}
    Proof. By using the invertibility of the matrix $G^{\tau}_N$, we define the controllers $v(i)$ and $u_1(i-\tau)$ as
    \begin{eqnarray}
        v(i)&=&D'C'(i-1)\cdots C'(0)(G^{\tau}_N)^{-1}x,\label{ct15}\\
        u_1(i-\tau)&=&D_1'E[C'(i-1)\cdots C'(0)|\mathcal{F}(i-\tau-1)]\nonumber\\
        &&\times (G^{\tau}_N)^{-1}x.\label{ct16}
    \end{eqnarray}
    It is easily verified that (\ref{ct11}) is solvable with $x(N+1)\in \mathcal{S}$ under the controllers (\ref{ct15})-(\ref{ct16}).
    Furthermore, from (\ref{ct12}) with $x(N+1)=0$, it yields that
    \begin{eqnarray}
        x(k)&=&E\Big[\sum_{i=k}^NC(k)C(k+1)\cdots C(i-1)\nonumber\\
        &&\times [Dv(i)+D_1u_1(k-\tau)]\Big|\mathcal{F}(k-1)\Big].\nonumber
    \end{eqnarray}
    This implies that
    \begin{eqnarray}
        x(0)\hspace{-2mm}&=&\hspace{-2mm}E\Big[\sum_{i=0}^NC(0)C(1)\cdots C(i-1)\big[DD'C'(i-1)\cdots\nonumber\\
        &&\hspace{-4mm}\times  C'(0)+D_1D_1'E[C'(i-1)\cdots C'(0)|\mathcal{F}(i-\tau-1)]\big]\Big]\nonumber\\
        &&\hspace{-4mm}\times(G^{\tau}_N)^{-1}x.\nonumber
    \end{eqnarray}
    Accordingly, under the controllers (\ref{ct15})-(\ref{ct16}), system (\ref{ct11}) satisfies $x(0)=x$ and $x(N+1)=0$, that is system (\ref{ct11}) is exactly null controllable. In addition, by applying Lemma (\ref{lem4}), system (\ref{ct11}) is also exactly controllable.
    The proof is now completed. \hfill $\blacksquare$

    \begin{remark}
        The sufficient condition obtained in Theorem \ref{thm4} contain that in Theorem \ref{lem3} as a special case. In fact, if $B_1=0$ in (\ref{ct10}), then the invertibility of the matrix Gramian matrix $G^{\tau}_N$ is reduced to that of that of $G_N$.
    \end{remark}


    \section{Stochastic Systems with State Delay}

    In this section, we consider the stochastic systems with state delay where the dynamics are governed by
    \begin{eqnarray}
        x(k+1)&=&[Ax(k)+A_1x(k-d)+Bu(k)]\nonumber\\
        &&\hspace{-3mm}+w(k)[\bar{A}x(k)+\bar{B}u(k)],\label{ct1}
    \end{eqnarray}
    where constant matrices $A, B, \bar{A}, \bar{B}$ have been defined in (\ref{c1}), $A_1$ is a constant matrices with compatible dimension, $d>0$ represents the state delay.

    By applying (\ref{a3}) and (\ref{c24}), (\ref{ct1}) becomes
    \begin{eqnarray}
        x(k+1)&=&[Ax(k)+A_1x(k-d)+Lq(k)+Fv(k)]\nonumber\\
        &&+w(k)[\bar{A}x(k)+q(k)]\nonumber\\
        &\triangleq&[(A-L\bar{A})x(k)+A_1x(k-d)+Lz(k)\nonumber\\
        &&+Fv(k)]+w(k)z(k)
    \end{eqnarray}
    where $z(k)=\bar{A}x(k)+q(k).$

    \begin{lemma}
        Under (\ref{a3}) and the invertibility of the matrix $A-L\bar{A}$, system (\ref{ct1}) is equivalently rewritten as
        \begin{eqnarray}
            x(k)&=&Cx(k+1)+C_1x(k-d)+\bar{C}z(k)\nonumber\\
            &&+Dv(k)-Cw(k)z(k),\label{ct2}
        \end{eqnarray}
        where $C, \bar{C}, D$ have the same definitions as those in Lemma \ref{lem1} and $C_1\triangleq -(A-L\bar{A})^{-1}A_1$.
        In particular, the solution to (\ref{ct2}) is given by
        \begin{eqnarray}
            x(k)&=&E\Big[C(k)x(k+1)+C_1x(k-d)\nonumber\\
            &&+Dv(k)\Big|\mathcal{F}(k-1)\Big],\label{ct3}
        \end{eqnarray}
        and $z(k)$ satisfies (\ref{c2}).
    \end{lemma}
    Proof. The proof follows similarly to those in Lemma \ref{lem1} and \ref{lem2}. So we omit the details. \hfill $\blacksquare$

    By introducing the backward iterations:
    \begin{eqnarray}
        P(k)&=&I, ~~k=N, \ldots, N-d+1,\nonumber\\
        P(k)&=&[I-CP(k+1)\cdots CP(k+d)C_1]^{-1},\nonumber\\
        &&k=N-d, \ldots, 0,\label{ct18}
    \end{eqnarray}
    we define the following Gramian matrix:
    \begin{eqnarray}
        G^d_{N}&=&E\Big[\sum_{j=0}^{N}P(0)C(0)\cdots P(j-1)C(j-1)\nonumber\\
        &&\times P(j)DD'P'(j)C'(j-1)P'(j-1)\cdots\nonumber\\
        &&\times C'(0)P'(0)\Big].\label{ct9}
    \end{eqnarray}

    We now give the controllability condition for system (\ref{ct1}).
    \begin{theorem}\label{thm5}
        Under (\ref{a3}), the invertibility of the matrix $A-L\bar{A}$ and the solvability of $P(k)$ in (\ref{ct18}), if there exists a positive integer $N$ such that the Gramian matrix $G^d_{N}$ in (\ref{ct9}) is invertible, then system (\ref{ct1}) with state delay is exactly null controllable, and is also exactly controllable, that is, for any nonzero $x\in R^l$ and $\xi^d\in \mathcal{S}^d$, there exist a positive integer $N$ and a sequence of controllers $\{u(k), k\in [0,N]\}$ which belongs to $l_{\mathcal{F}}^2([0, N], R^m)$ such that $x(0)=x$ and $x(N+1)=\xi^d\in \mathcal{S}^d$ where \begin{eqnarray}
            \mathcal{S}^d&=&\{\xi^d|\xi^d\in L^2(\Omega, R^s)~ such~that~BSDE\nonumber\\
            &&x^d(k)=Cx^d(k+1)+C_1x(k-d)+\bar{C}z^d(k)\nonumber\\
            &&-C_1w(k)z^d(k)~with~x^d(N+1)=\xi^d\nonumber\\
            &&is~solvable\}.\nonumber
        \end{eqnarray}
    \end{theorem}
    Proof. By using the invertibility of the Gramian matrix $G^d_{N}$, we define the controller:
    \begin{eqnarray}
        v(i)&=&D'P'(j)C'(j-1)P'(j-1)\cdots\nonumber\\
        &&\times C'(0)P'(0)(G^d_N)^{-1}x.\label{ct17}
    \end{eqnarray}
    It can be verified that (\ref{ct2}) is solvable under the controller (\ref{ct17}). Furthermore, from (\ref{ct3}) with $x(N+1)=0$, it yields for $k=N,\cdots, N-d+1$ that
    \begin{eqnarray}
        x(k)&=&E\Big[\sum_{j=k}^{N}P(k)C\cdots P(j-1)CP(j)C_1x(j-d)\nonumber\\
        &&+\sum_{j=k}^{N}P(k)C(k)\cdots P(j-1)C(j-1)\nonumber\\
        &&\times P(j)Dv(j)\Big|\mathcal{F}(k-1)\Big].\label{ct4}
    \end{eqnarray}
    Taking continuous iterations yields for $k=N-d, \cdots, 0$ that
    \begin{eqnarray}
        x(k)&=&E\Big[\sum_{j=k}^{k+d-1}P(k)C\cdots P(j-1)CP(j)C_1x(j-d)\nonumber\\
        &&+\sum_{j=k}^{N}P(k)C(k)\cdots P(j-1)C(j-1)\nonumber\\
        &&\times P(j)Dv(j)\Big|\mathcal{F}(k-1)\Big].\label{ct5}
    \end{eqnarray}
    Thus, by using the controller (\ref{ct17}) and $x(s)=0$ for $s=-1, \ldots, -d$, it follows from (\ref{ct5}) that
    \begin{eqnarray}
        x(0)&=&E\Big[\sum_{j=0}^{N}P(0)C(0)\cdots P(j-1)C(j-1)\nonumber\\
        &&\times P(j)Dv(j)\Big]\nonumber\\
        &=&E\Big[\sum_{j=0}^{N}P(0)C(0)\cdots P(j-1)C(j-1)\nonumber\\
        &&\times P(j)DD'P'(j)C'(j-1)P'(j-1)\cdots\nonumber\\
        &&\times C'(0)P'(0)\Big](G^d_N)^{-1}x.\nonumber
    \end{eqnarray}
    That is, system (\ref{ct1}) satisfies $x(0)=x$ and $x(N+1)=0$ under the controller (\ref{ct17}).
    In addition, by applying similar proof to that in Lemma (\ref{lem4}), system (\ref{ct1}) is also exactly controllable.
    The proof is now completed. \hfill $\blacksquare$

    \begin{remark}
        The sufficient condition obtained in Theorem \ref{thm5} contain that in Theorem \ref{lem3} as a special case. In fact, if $A_1=0$ in (\ref{ct10}), then $C_1=0$ and $P(k)=I$ for $k=N, \ldots, 0$ which implies that $G^d_N=G_N$, that is, the invertibility of the matrix Gramian matrix $G^{d}_N$ is reduced to that of that of $G_N$.
    \end{remark}

    \section{Numerical Examples}
In this section, we present four numerical examples to illustrate the effectiveness of the proposed results.
\begin{example}\label{ex1}
 Consider system (\ref{c1}) with parameters given by
   \begin{eqnarray}
    &&A=\left[
    \begin{array}{cc}
        1 \ & \ 1 \\
        -1 \ & \ -2 \\
    \end{array}
    \right], \bar{A}=\left[
    \begin{array}{cc}
        1 \ & \ 2 \\
        0 \ & \ 1 \\
    \end{array}
    \right],
        \nonumber\\
&&B=\left[
\begin{array}{ccc}
    1 \ & \ 2  \ & \ -2\\
    1 \ & \ 2  \ & \ -3\\
\end{array}
\right],\bar{B}=\left[
    \begin{array}{ccc}
        1 \ & \ 1  \ & \ 0\\
        0 \ & \ 1  \ & \ -2\\
    \end{array}
    \right], \nonumber\\
    &&n=2,m=3,N=2.\nonumber
\end{eqnarray}

In this case, $M=\left[
\begin{array}{ccc}
    1 \ & \ -1  \ & \ 2\\
    0 \ & \ 1  \ & \ -2\\
    0 \ & \ 0  \ & \ -1\\
\end{array}
\right]$. From (\ref{c9}) and (\ref{c48}), it can be verified by direct calculation that $Rank(G_N)=2$ and $Rank(R)=2$. By applying  Theorem \ref{thm2}, we have that the stochastic system (\ref{c1}) is exactly controllable.
\end{example}

\begin{example}\label{ex2}
    Consider system (\ref{c1}) with parameters given by
    \begin{eqnarray}
        &&A=\left[
        \begin{array}{cc}
            1 \ & \ 2 \\
            2 \ & \ 1 \\
        \end{array}
        \right], \bar{A}=\left[
        \begin{array}{cc}
            1 \ & \ 0 \\
            0 \ & \ 1 \\
        \end{array}
        \right], \nonumber\\
        &&B=\left[
        \begin{array}{ccc}
            1 \ & \ 2  \ & \ -1\\
            1 \ & \ 2  \ & \ 0\\
        \end{array}
        \right], \bar{B}=\left[
        \begin{array}{ccc}
            1 \ & \ 1  \ & \ 0\\
            0 \ & \ 1  \ & \ -1\\
        \end{array}
        \right], \nonumber\\
        &&n=2,m=3,l=1,N=2.\nonumber
    \end{eqnarray}

    In this case, $M=\left[
    \begin{array}{ccc}
        1 \ & \ -1  \ & \ -1\\
        0 \ & \ 1  \ & \ 1\\
        0 \ & \ 0  \ & \ 1\\
    \end{array}
    \right]$ and $H=\left[
    \begin{array}{cc}
        1 & 1
    \end{array}
    \right]$. By solving $\bar{G}_N$ and $\bar{R}$, it is calculated that $\bar{G}_N=31$ and $\bar{R}=\left[
    \begin{array}{c}
        -1  \ \ \  -1  \ \ \  2  \ \ \  -1   \ \ \   -4   \ \ \    2  \ \ \   2
    \end{array}
    \right]$. This further gives that
    $Rank(\bar{G}_N)=1$ and $Rank(\bar{R})=1$.
 Therefore, by Theorem \ref{thm3}, the stochastic system (\ref{c1}) is $H$-partially exactly controllable.
\end{example}

\begin{example}\label{ex3}
    Consider system (\ref{ct10}) with parameters given by
    \begin{eqnarray}
        &&A=\left[
        \begin{array}{cc}
            1 \ & \ -1 \\
            0 \ & \ 1 \\
        \end{array}
        \right],  \bar{A}=\left[
        \begin{array}{cc}
            1 \ & \ 2 \\
            1 \ & \ 1 \\
        \end{array}
        \right],  \nonumber\\
        && B=\left[
        \begin{array}{ccc}
            1 \ & \ 0  \ & \ 1\\
            1 \ & \ 1  \ & \ 0\\
        \end{array}
        \right], \bar{B}=\left[
        \begin{array}{ccc}
            1 \ & \ 0  \ & \ 0\\
            0 \ & \ 1  \ & \ 0\\
        \end{array}
        \right], \nonumber\\
        &&{B}_1=\left[
        \begin{array}{ccc}
            0 \ & \ 1  \ & \ 0\\
            1 \ & \ 0  \ & \ 1\\
        \end{array}
        \right],n=2,m=3,N=2, \tau=1.\nonumber
    \end{eqnarray}

In this case, $M=I$. By using Theorem \ref{thm4}, it is calculated from (\ref{ct14}) that $Rank(G^{\tau}_N)=2$, that is, the Gramian matrix $G^{\tau}_{N}$ is invertible. This indicates that system (\ref{ct10}) with input delay is exactly controllable.
\end{example}

\begin{example}\label{ex4}
    Consider system (\ref{ct1}) with parameters given by
    \begin{eqnarray}
        &&A=\left[
        \begin{array}{cc}
            2 \ & \ 0 \\
            0 \ & \ -1 \\
        \end{array}
        \right], \bar{A}=\left[
        \begin{array}{cc}
            1 \ & \ 2 \\
            0 \ & \ 2 \\
        \end{array}
        \right], \nonumber\\
        && B=\left[
        \begin{array}{ccc}
            1 \ & \ -1  \ & \ 2\\
            1 \ & \ 1  \ & \ 0\\
        \end{array}
        \right], \bar{B}=\left[
        \begin{array}{ccc}
            1 \ & \ 0  \ & \ 0\\
            0 \ & \ 1  \ & \ 0\\
        \end{array}
        \right], \nonumber\\
        &&A_1=\left[
        \begin{array}{cc}
            1 \ & \ 0 \\
            0  \ & \ 1\\
        \end{array}
        \right],n=2,m=3,N=2, d=1.\nonumber
    \end{eqnarray}

    In this case, $M=I$. By applying Theorem \ref{thm5}, it yields from (\ref{ct9}) that $Rank(G^{d}_N)=2$, that is, the Gramian matrix $G^{d}_{N}$ is invertible. Thus, the system (\ref{ct1}) with state delay is exactly controllable.
\end{example}

    \section{Conclusions}

    In this paper, we studied the exact controllability of discrete-time stochastic system.
    To overcome the adaptiveness constraint of the controllers and the solvability challenging of stochastic difference equation with terminal value, we innovatively transformed the FSDE into BSDE. Accordingly, we obtained both the Gramian matrix criterion and the Rank criterion
    for the exact controllability of discrete-time system with multiplicative noise. Furthermore, we derived the controllability conditions for the cases with output measurements, input delay and state delay.


\end{document}